\documentclass[12pt,oneside,english]{amsart}
\usepackage[T1]{fontenc}
\usepackage[latin9]{inputenc}
\usepackage{babel}
\usepackage{textcomp}
\usepackage{mathtools}
\usepackage{amstext}
\usepackage{amsthm}
\usepackage{amssymb}
\usepackage[unicode=true,pdfusetitle,
 bookmarks=true,bookmarksnumbered=false,bookmarksopen=false,
 breaklinks=false,pdfborder={0 0 1},backref=false,colorlinks=false]
 {hyperref}

\makeatletter
\numberwithin{equation}{section}
\theoremstyle{plain}
\newtheorem{conjecture}{\protect\conjecturename}[section]
\theoremstyle{plain}
\newtheorem{thm}{\protect\theoremname}[section]
\theoremstyle{plain}
\newtheorem{lem}{\protect\lemmaname}[section]
\theoremstyle{plain}
\newtheorem{prop}{\protect\propositionname}[section]

\usepackage{booktabs}
\usepackage{textcomp}
\usepackage{listings}
\usepackage{xcolor}

\lstdefinelanguage{Maple}{
  morekeywords={proc,local,for,to,do,end,while,if,then,elif,else,error,
    printf,print,restart,seq,add,expand,diff,degree,coeff,eval,subs,min,max,
    floor,ceil,sqrt,powmod,numelems,Array,Vector,Matrix,table,type,
    assigned,seq,return,mod},
  sensitive=true,
  morecomment=[l]{\#},
  morestring=[b]",
  alsoletter={_}
}

\makeatother

\providecommand{\conjecturename}{Conjecture}
\providecommand{\lemmaname}{Lemma}
\providecommand{\propositionname}{Proposition}
\providecommand{\theoremname}{Theorem}

\begin{document}
\title[On the Baruah\textendash Sarma Congruence Conjecture]{A Proof of the Baruah\textendash Sarma Congruence Conjecture for the
Rogers\textendash Ramanujan Continued Fraction}
\author{Bing He}
\thanks{School of Mathematics and Statistics, HNP-LAMA, Central South University
\\
Changsha 410083, Hunan, People's Republic of China}
\thanks{yuhelingyun@foxmail.com; yuhe123456@foxmail.com}
\begin{abstract}
  Recently, Baruah and Sarma studied congruences satisfied by four
infinite products closely related to the celebrated Rogers\textendash Ramanujan
continued fraction and proposed a conjecture consisting of three families
of congruences. In this paper, we prove this conjecture. Our proof
is modular: we realize the relevant quotients as modular functions,
construct associated holomorphic modular forms, and apply a Sturm-type
theorem. This approach differs from the methods used by Baruah and
Sarma in their study of these products.
\end{abstract}

\maketitle

\section{Introduction}

 The celebrated Rogers\textendash Ramanujan continued fraction is
defined by

\[
r(\tau):=\cfrac{q^{1/5}}{1+\cfrac{q}{1+\cfrac{q^{2}}{1+\cfrac{q^{3}}{1+\ddots}}}},\quad q=e^{2\pi i\tau},\;\Im\tau>0.
\]
It was first discovered by Rogers \cite{R1} in 1894 and rediscovered
independently by Ramanujan between 1910 and 1913. The product representation
\[
r(\tau)=q^{1/5}\frac{(q;q^{5})_{\infty}(q^{4};q^{5})_{\infty}}{(q^{2};q^{5})_{\infty}(q^{3};q^{5})_{\infty}},\qquad(z;q)_{\infty}\coloneqq\prod_{k=0}^{\infty}(1-zq^{k}),
\]
 first proved by Rogers \cite{R1} and later by Ramanujan \cite{R},
shows that $r(\tau)$ is essentially a quotient of the two Rogers\textendash Ramanujan
functions 
\[
G(q)=\sum_{n\geq0}\frac{q^{n^{2}}}{(q;q)_{n}}=\frac{1}{(q;q^{5})_{\infty}(q^{4};q^{5})_{\infty}},
\]
and
\[
H(q)=\sum_{n\geq0}\frac{q^{n^{2}+n}}{(q;q)_{n}}=\frac{1}{(q^{2};q^{5})_{\infty}(q^{3};q^{5})_{\infty}},
\]
where
\[
(z;q)_{n}=\frac{(z;q)_{\infty}}{(zq^{n};q)_{\infty}}.
\]
 Combinatorially, $G(q)$ and $H(q)$ are the generating functions
for partitions into parts congruent to $\pm1(\bmod\:5)$ and $\pm2(\bmod\:5),$
respectively. Thus the coefficients of $r(\tau)$ and of its reciprocal
are natural partition-theoretic arithmetic objects. Andrews \cite{A}
provided partition-theoretic interpretations for these coefficients,
while Richmond\textendash Szekeres \cite{RS} and Hirschhorn \cite{H}
studied their sign patterns and periodicity.

Recently, Baruah and Sarma \cite{BS} shifted attention from $r(\tau)$
itself to fifth powers of the Rogers\textendash Ramanujan quotient
\[
\mathcal{R}(q)=\frac{(q,q^{4};q^{5})_{\infty}}{(q^{2},q^{3};q^{5})_{\infty}},
\]
so that $r(\tau)=q^{1/5}\mathcal{R}(q).$ They studied the coefficient
sequences of 
\[
\frac{1}{\mathcal{R}^{5}(q)},\quad\ensuremath{\mathcal{R}^{5}(q)},\quad\ensuremath{\frac{\mathcal{R}^{5}(q)}{\mathcal{R}(q^{5})}},\quad\frac{\mathcal{R}(q^{5})}{\mathcal{R}^{5}(q)},
\]
proved periodic sign patterns, and discovered several new congruences.
In particular, writing

\[
C(q):=\frac{\mathcal{R}^{5}(q)}{\mathcal{R}(q^{5})}=:\sum_{n\geq0}C(n)q^{n},\qquad D(q):=\frac{\mathcal{R}(q^{5})}{\mathcal{R}^{5}(q)}=:\sum_{n\geq0}D(n)q^{n},
\]
Baruah and Sarma proposed the following conjecture based on numerical
calculations.
\begin{conjecture}
\label{conj1} For every integer $n\ge0$, we have

\[
C(27n+18)\equiv D(27n+18)\equiv0(\bmod\:3),
\]
\[
C(16n+12)\equiv D(16n+12)\equiv0(\bmod\:4),
\]
\[
C(32n+28)\equiv D(32n+28)\equiv0(\bmod\:8).
\]
\end{conjecture}
  From a combinatorial perspective, $C(q)$ and $D(q)$ are not
ad hoc $q$-series: they are quotients built from the Rogers\textendash Ramanujan
partition generating functions $G(q)$ and $H(q),$ together with
their fifth powers and the $q^{5}$-substitution. Consequently, congruences
for $C(n)$ and $D(n)$ are arithmetic statements about partition-theoretic
quantities, not merely formal identities for infinite products. The
conjecture is also a natural next layer in the study of the Rogers\textendash Ramanujan
quotient: earlier work concerned signs and periodicity for $r(\tau)$
and its reciprocal, whereas the fifth-power quotients above involve
a more entangled product structure whose arithmetic is not implied
by the earlier sign-pattern results. Proving Conjecture \ref{conj1}
therefore completes an important part of the arithmetic picture for
fifth powers of $\mathcal{R}(q),$ and it does so for three moduli
3, 4, 8 simultaneously.

To the best of our knowledge, no proof of Conjecture \ref{conj1}
has been given so far.

In this paper, we supply such a proof by applying the theory of modular
forms.
\begin{thm}
\label{thm:main} Conjecture \ref{conj1} is true.
\end{thm}
The proof is modular and uniform in the three congruence families.
Define 
\[
\mathcal{C}(\tau)=\frac{r(\tau)^{5}}{r(5\tau)},\qquad\mathcal{D}(\tau)=\frac{r(5\tau)}{r(\tau)^{5}}.
\]
In Section \ref{sec:auxiliary}, we show that $\mathcal{C}$ and $\mathcal{D}$
are modular functions on a congruence subgroup $G$ of level 25, with
neither zeros nor poles in the upper half-plane, and we determine
their orders at the cusps. For 
\[
m\in\{27,16,32\},\qquad t=\left\lceil \frac{m}{5}\right\rceil ,\qquad X\in\{\mathcal{C},\mathcal{D}\},
\]
we multiply $X$ by a suitable power of the discriminant, 
\[
H_{X}(\tau)=\Delta(5m\tau)^{t}X(\tau),\qquad\Delta(\tau)=\eta(\tau)^{24},
\]
and prove that $H_{X}$ is a holomorphic modular form of weight $12t$
on $\ensuremath{G\cap\Gamma_{0}(25m)}.$ Applying a projection operator
that selects coefficients in a fixed residue class modulo $\ensuremath{m},$
we obtain a holomorphic modular form on 
\[
G_{m}=G\cap\Gamma_{0}(25m^{2})
\]
with integer Fourier coefficients. In Section \ref{sec:proof}, we
apply a Sturm-type theorem for congruence subgroups to reduce the
desired congruences to a finite verification. More precisely, for
each triple 
\[
(m,s,\ell)\in\{(27,18,3),(16,12,4),(32,28,8)\},
\]
 it suffices to check the coefficients of the projected form through
an explicit Sturm bound $\ensuremath{B=t\mu_{m}},$ where $\ensuremath{\mu_{m}=[\mathrm{SL}_{2}(\mathbb{Z}):G_{m}]}.$
These finite checks are carried out exactly modulo $24$ and are recorded
in Table \ref{tab:cert}; the corresponding Maple program is given
in Appendix \ref{sec:Verification}. This proves Theorem \ref{thm:main}.

The paper is organized as follows. Section \ref{sec:auxiliary} contains
the auxiliary modularity and cusp-order results, together with the
construction and projection of the forms $\ensuremath{H_{X}}.$ Section
\ref{sec:proof} gives the proof of Theorem \ref{thm:main} via a
Sturm-type theorem and the finite verifications. Appendix \ref{sec:Verification}
provides the verification program used for Table \ref{tab:cert}.

\section{\label{sec:auxiliary} Auxiliary results}

The coefficients $C(n)$ and $D(n)$ in Conjecture \ref{conj1} are
not arbitrary $q$-series coefficients: they arise from quotients
of the Rogers\textendash Ramanujan partition generating functions
$G(q)$ and $H(q),$ together with their fifth powers and the $q^{5}$-substitution.
Congruences for such coefficients cannot be obtained by purely formal
$q$-series manipulation. The purpose of this section is to develop
a modular-form framework which converts the discrete combinatorial
congruences into statements about modular functions, their cusp orders,
and projections of modular forms. The key point is that the modular
functions $\mathcal{C}$ and $\mathcal{D}$ below have very mild pole
behaviour at the cusps, so that a single discriminant factor $\Delta(5m\tau)^{t}$
compensates all possible poles uniformly for the three moduli 3, 4,
8. This is what makes a uniform modular proof of all three families
in Conjecture \ref{conj1} possible.

For an integer $\ensuremath{N\ge1},$ define
\[
\Gamma(N)=\left\{ \begin{pmatrix}a & b\\
c & d
\end{pmatrix}\in\mathrm{SL}_{2}(\mathbb{Z}):a\equiv d\equiv1(\bmod\:N),\;b\equiv c\equiv0(\bmod\:N)\right\} ,
\]
\[
\Gamma_{1}(N)=\left\{ \begin{pmatrix}a & b\\
c & d
\end{pmatrix}\in\mathrm{SL}_{2}(\mathbb{Z}):a\equiv d\equiv1(\bmod\:N),\quad c\equiv0(\bmod\:N)\right\} ,
\]
\[
\Gamma_{0}(N)=\left\{ \begin{pmatrix}a & b\\
c & d
\end{pmatrix}\in\mathrm{SL}_{2}(\mathbb{Z}):c\equiv0(\bmod\:N)\right\} ,
\]
and
\[
G=\left\{ \begin{pmatrix}a & b\\
c & d
\end{pmatrix}\in\mathrm{SL}_{2}(\mathbb{Z}):25|c,\;a\equiv d\equiv1(\bmod\:5)\;\mathrm{or}\;a\equiv d\equiv-1(\bmod\:5)\right\} .
\]

\begin{lem}
\label{lem:2-1} The functions 
\[
\mathcal{C}(\tau)=\frac{r(\tau)^{5}}{r(5\tau)},\qquad\mathcal{D}(\tau)=\frac{r(5\tau)}{r(\tau)^{5}}
\]
 are modular functions on $G$. They have neither zeros nor poles
in the upper half-plane.
\end{lem}
\begin{proof}
The classical transformation laws \cite[Proposition 2]{D} for the
Rogers\\
\textendash Ramanujan function give 
\[
r(\gamma\tau)=r(\tau),\quad\gamma\in\Gamma(5),
\]
and 
\[
r(\tau+1)=\zeta_{5}r(\tau),\qquad\zeta_{5}=e^{2\pi i/5}.
\]
Therefore 
\[
r(\tau+1)^{5}=r(\tau)^{5}.
\]
Since $\Gamma_{1}(5)/\Gamma(5)$ is generated by the class of the
translation matrix, it follows that 
\begin{equation}
r(\gamma\tau)^{5}=r(\tau)^{5}\qquad\text{for all }\gamma\in\Gamma_{1}(5).\label{eq:1}
\end{equation}

Now let 
\[
\gamma=\begin{pmatrix}a & b\\
c & d
\end{pmatrix}\in G.
\]
 Replacing $\gamma$ by $-\gamma$ if necessary, we may assume 
\[
a\equiv d\equiv1(\bmod\:5).
\]
 This replacement does not change the action on the upper half-plane,
since $-I$ acts trivially on modular functions of weight $\ensuremath{0}.$
Because $\ensuremath{25|c},$ we have $\ensuremath{c\equiv0(\bmod\:5)},$
and $\ensuremath{\gamma\in\Gamma_{1}(5)}.$ By \eqref{eq:1}, we have
\begin{equation}
r(\gamma\tau)^{5}=r(\tau)^{5}.\label{eq:2}
\end{equation}
Next observe that 
\[
5\gamma\tau=\frac{a(5\tau)+5b}{(c/5)(5\tau)+d}.
\]
 Consider the matrix 
\[
M=\begin{pmatrix}a & 5b\\
c/5 & d
\end{pmatrix}.
\]
 Since $\ensuremath{25|c},$ the entry $c/5$ is an integer, so $M$
is integral. Moreover, 
\[
\det M=ad-5b\cdot\frac{c}{5}=ad-bc=1,
\]
so $\ensuremath{M\in\mathrm{SL}_{2}(\mathbb{Z})}.$ Modulo $\ensuremath{5},$
\[
a\equiv1,\qquad d\equiv1,\qquad5b\equiv0,\qquad\frac{c}{5}\equiv0(\bmod\:5),
\]
 because $\ensuremath{25|c}.$ Hence $\ensuremath{M\in\Gamma(5)}.$
Therefore, using $r(\gamma\tau)=r(\tau)$ for $\ensuremath{\gamma\in\Gamma(5)},$
we obtain 
\begin{equation}
r(5\gamma\tau)=r(M(5\tau))=r(5\tau).\label{eq:3}
\end{equation}

Combining \eqref{eq:2} and \eqref{eq:3}, we obtain 
\[
\mathcal{C}(\gamma\tau)=\frac{r(\gamma\tau)^{5}}{r(5\gamma\tau)}=\frac{r(\tau)^{5}}{r(5\tau)}=\mathcal{C}(\tau),
\]
 and 
\[
\mathcal{D}(\gamma\tau)=\frac{r(5\gamma\tau)}{r(\gamma\tau)^{5}}=\frac{r(5\tau)}{r(\tau)^{5}}=\mathcal{D}(\tau).
\]
 Thus $\mathcal{C}$ and $\mathcal{D}$ are invariant under $\ensuremath{G}.$

It remains to check that they have neither zeros nor poles in the
upper half-plane. For $\ensuremath{|q|<1},$ every factor $1-q^{m}$
appearing in the product representation of $\mathcal{R}(q)$ is nonzero,
hence $\ensuremath{\mathcal{R}(q)\neq0}.$ Since $q\neq0$ in the
upper half-plane, we have $\ensuremath{r(\tau)=q^{1/5}\mathcal{R}(q)\neq0},$
and similarly $\ensuremath{r(5\tau)\neq0}.$ Therefore, both numerator
and denominator in 
\[
\mathcal{C}(\tau)=\frac{r(\tau)^{5}}{r(5\tau)},\qquad\mathcal{D}(\tau)=\frac{r(5\tau)}{r(\tau)^{5}}
\]
are holomorphic and nonzero in the upper half-plane. Consequently
$\mathcal{C}$ and $\mathcal{D}$ are holomorphic and nonzero there;
in particular, they have neither zeros nor poles.

Finally, $G$ is a discrete subgroup of $\ensuremath{\mathrm{SL}_{2}(\mathbb{Z})},$
and $\mathcal{C},\mathcal{D}$ are meromorphic functions on the upper
half-plane invariant under $\ensuremath{G}.$ Hence they are modular
functions on $\ensuremath{G}.$ This proves Lemma\ref{lem:2-1}.
\end{proof}
Combinatorially, Lemma \ref{lem:2-1} says that the arithmetic of
the partition-theoretic coefficient sequences $C(n)$ and $D(n)$
can be controlled by modular functions on the congruence group $G.$
The next step is to examine how mild their singularities are at the
cusps; this is precisely what determines whether a single discriminant
factor can compensate them.

For a form $f$ of weight $k$ and $\sigma\in\mathrm{SL}_{2}(\mathbb{Z})$
taking $\infty$ to $a/c$, let $v_{a/c}(f)$ denote the leading exponent
in the $q$-expansion of $f|_{k}\sigma$. This exponent has not been
multiplied by the cusp width; multiplication by the positive cusp
width gives the usual local order, so its sign determines holomorphy.

Define 
\[
\eta_{N,g}(\tau)=q^{\frac{N}{2}B_{2}(g/N)}(q^{g},q^{N-g};q^{N})_{\infty},\qquad B_{2}(x)=x^{2}-x+\frac{1}{6},
\]
and put $P_{2}(x)=B_{2}(\{x\})$, where $\{x\}$ is the fractional
part of $x$. The generalized eta cusp formula \cite[Lemma 2.9]{G}
is 
\begin{equation}
v_{a/c}(\eta_{N,g})=\frac{\gcd(c,N)^{2}}{2N}P_{2}\!\left(\frac{ag}{\gcd(c,N)}\right).\label{eq:eta-cusp}
\end{equation}
 We use $\gcd(0,N)=N$ at $\infty$. Here 
\[
r(\tau)=\frac{\eta_{5,1}}{\eta_{5,2}},\qquad r(5\tau)=\frac{\eta_{25,5}}{\eta_{25,10}}.
\]
 For $5\nmid a$, put 
\[
\epsilon(a)=\left(\frac{a}{5}\right),
\]
 where $\left(\frac{\cdot}{5}\right)$ denotes the Legendre symbol
modulo $5$.
\begin{lem}
\label{lem:2-2} For $\ensuremath{5\nmid a},$ we have 
\[
v_{a/c}(\mathcal{C})=\begin{cases}
\epsilon(a), & 5|c\ \text{and}\ 25\nmid c,\\
0, & \text{otherwise},
\end{cases}\qquad v_{a/c}(\mathcal{D})=-v_{a/c}(\mathcal{C}).
\]
 Consequently, a negative exponent can occur only when $5|c$ and
$\ensuremath{25\nmid c},$ and it is always $\ensuremath{-1}.$
\end{lem}
\begin{proof}
Since 
\[
\mathcal{C}(\tau)=\frac{r(\tau)^{5}}{r(5\tau)}=\frac{(\eta_{5,1}/\eta_{5,2})^{5}}{\eta_{25,5}/\eta_{25,10}},
\]
 we have 
\begin{equation}
v_{a/c}(\mathcal{C})=5\bigl(v_{a/c}(\eta_{5,1})-v_{a/c}(\eta_{5,2})\bigr)-\bigl(v_{a/c}(\eta_{25,5})-v_{a/c}(\eta_{25,10})\bigr).\label{eq:2-1}
\end{equation}
 We compute the four terms using \eqref{eq:eta-cusp}. First,
\[
v_{a/c}(\eta_{5,1})=\frac{\gcd(c,5)^{2}}{10}P_{2}\!\left(\frac{a}{\gcd(c,5)}\right),
\]
and 
\[
v_{a/c}(\eta_{5,2})=\frac{\gcd(c,5)^{2}}{10}P_{2}\!\left(\frac{2a}{\gcd(c,5)}\right).
\]
 Hence 
\begin{equation}
\begin{aligned} & v_{a/c}(\eta_{5,1})-v_{a/c}(\eta_{5,2})\\
 & =\frac{\gcd(c,5)^{2}}{10}\left[P_{2}\!\left(\frac{a}{\gcd(c,5)}\right)-P_{2}\!\left(\frac{2a}{\gcd(c,5)}\right)\right].
\end{aligned}
\label{eq:2-2}
\end{equation}

For $\ensuremath{N=25},$ 
\[
v_{a/c}(\eta_{25,5})=\frac{\gcd(c,25)^{2}}{50}P_{2}\!\left(\frac{5a}{\gcd(c,25)}\right),
\]
and 
\[
v_{a/c}(\eta_{25,10})=\frac{\gcd(c,25)^{2}}{50}P_{2}\!\left(\frac{10a}{\gcd(c,25)}\right).
\]
 Hence 
\begin{equation}
\begin{aligned} & v_{a/c}(\eta_{25,5})-v_{a/c}(\eta_{25,10})\\
 & =\frac{\gcd(c,25)^{2}}{50}\left[P_{2}\!\left(\frac{5a}{\gcd(c,25)}\right)-P_{2}\!\left(\frac{10a}{\gcd(c,25)}\right)\right].
\end{aligned}
\label{eq:2-3}
\end{equation}

We now distinguish three cases according to the power of $5$ dividing
$\ensuremath{c}.$

\medskip \noindent \textbf{Case 1: \(5\nmid c\).} Then $\gcd(c,5)=1$
and $\ensuremath{\gcd(c,25)=1}.$ From \eqref{eq:2-2}, we get 
\[
v_{a/c}(\eta_{5,1})-v_{a/c}(\eta_{5,2})=\frac{1}{10}\bigl(P_{2}(a)-P_{2}(2a)\bigr).
\]
 Since $a$ is an integer, we have $\ensuremath{P_{2}(a)=P_{2}(2a)=B_{2}(0)=1/6},$
so this difference is $\ensuremath{0}.$ From \eqref{eq:2-3} we derive
that 
\[
v_{a/c}(\eta_{25,5})-v_{a/c}(\eta_{25,10})=\frac{1}{50}\bigl(P_{2}(5a)-P_{2}(10a)\bigr)=0,
\]
 because $5a$ and $10a$ are integers. Substituting these into \eqref{eq:2-1}
gives
\[
v_{a/c}(\mathcal{C})=5\cdot0-0=0.
\]

\medskip \noindent \textbf{Case 2: \(5\mid c\) and \(25\nmid c\).}
Then $\gcd(c,5)=5$ and $\ensuremath{\gcd(c,25)=5}.$ From \eqref{eq:2-2}
we have 
\[
v_{a/c}(\eta_{5,1})-v_{a/c}(\eta_{5,2})=\frac{25}{10}\left[P_{2}\!\left(\frac{a}{5}\right)-P_{2}\!\left(\frac{2a}{5}\right)\right].
\]
 From \eqref{eq:2-3} we deduce that 
\[
v_{a/c}(\eta_{25,5})-v_{a/c}(\eta_{25,10})=\frac{25}{50}\left[P_{2}(a)-P_{2}(2a)\right]=0.
\]
 Substituting these into \eqref{eq:2-1} yields 
\begin{equation}
v_{a/c}(\mathcal{C})=\frac{25}{2}\left[P_{2}\!\left(\frac{a}{5}\right)-P_{2}\!\left(\frac{2a}{5}\right)\right].\label{eq:2-4}
\end{equation}
 We use the elementary identity 
\begin{equation}
P_{2}\!\left(\frac{a}{5}\right)-P_{2}\!\left(\frac{2a}{5}\right)=\frac{2\epsilon(a)}{25}\qquad(5\nmid a).\label{eq:P2-identity}
\end{equation}
Substituting \eqref{eq:P2-identity} into \eqref{eq:2-4} gives 
\[
v_{a/c}(\mathcal{C})=\frac{25}{2}\cdot\frac{2\epsilon(a)}{25}=\epsilon(a).
\]

\medskip \noindent \textbf{Case 3: \(25\mid c\).} Then $\gcd(c,5)=5$
and $\ensuremath{\gcd(c,25)=25}.$ From \eqref{eq:2-2} we get 
\[
v_{a/c}(\eta_{5,1})-v_{a/c}(\eta_{5,2})=\frac{25}{10}\left[P_{2}\!\left(\frac{a}{5}\right)-P_{2}\!\left(\frac{2a}{5}\right)\right].
\]
 From \eqref{eq:2-3} we obtain 
\begin{align*}
v_{a/c}(\eta_{25,5})-v_{a/c}(\eta_{25,10}) & =\frac{625}{50}\left[P_{2}\!\left(\frac{a}{5}\right)-P_{2}\!\left(\frac{2a}{5}\right)\right]\\
 & =\frac{25}{2}\left[P_{2}\!\left(\frac{a}{5}\right)-P_{2}\!\left(\frac{2a}{5}\right)\right].
\end{align*}
 Substituting these into \eqref{eq:2-1} produces 
\[
v_{a/c}(\mathcal{C})=5\cdot\frac{25}{10}\left[P_{2}\!\left(\frac{a}{5}\right)-P_{2}\!\left(\frac{2a}{5}\right)\right]-\frac{25}{2}\left[P_{2}\!\left(\frac{a}{5}\right)-P_{2}\!\left(\frac{2a}{5}\right)\right]=0.
\]
This proves the stated formula for $\mathcal{C}.$ Finally, since
\[
\mathcal{D}(\tau)=\frac{r(5\tau)}{r(\tau)^{5}}=\frac{1}{\mathcal{C}(\tau)},
\]
 we have 
\[
v_{a/c}(\mathcal{D})=-v_{a/c}(\mathcal{C}).
\]
 The last assertion follows because $\ensuremath{\epsilon(a)\in\{\pm1\}}.$
This proves Lemma \ref{lem:2-2}.
\end{proof}
Lemma \ref{lem:2-2} is the technical heart of the uniform argument:
away from the cusps with $5|c$ and $25\nmid c,$ the functions $\mathcal{C}$
and $\mathcal{D}$ are already holomorphic; at the exceptional cusps
their order is only $-1.$ Thus all possible poles are extremely mild
and can be compensated by a single discriminant factor. This is the
modular explanation for why the three congruence families in Conjecture
\ref{conj1} can be treated simultaneously.

Let $m\in\{27,16,32\}$ and set 
\[
t=\left\lceil \frac{m}{5}\right\rceil ,\qquad\Delta(\tau)=\eta(\tau)^{24}=q\prod_{j\ge1}(1-q^{j})^{24}.
\]
 For $X=\mathcal{C}$ or $X=\mathcal{D}$, define 
\begin{equation}
H_{X}(\tau)=\Delta(5m\tau)^{t}X(\tau).\label{eq:HX}
\end{equation}
The factor $\Delta(5m\tau)^{t}$ is chosen precisely to cancel the
possible order $-1$ poles described in Lemma \ref{lem:2-2} while
preserving the congruence-class structure modulo $m;$ its exponents
are multiples of $5m,$ which will be compatible with the projection
operator introduced below.
\begin{lem}
\label{lem:2-3} For each $m\in\{27,16,32\}$ and each $\ensuremath{X\in\{\mathcal{C},\mathcal{D}\}},$
the function $H_{X}$ is a holomorphic modular form of weight $12t$
on $\ensuremath{G\cap\Gamma_{0}(25m)}.$
\end{lem}
\begin{proof}
We split the proof into three parts: the transformation law, holomorphy
in the upper half-plane, and holomorphy at the cusps.

\medskip \noindent \textbf{(1) Transformation law.} The discriminant
$\Delta(\tau)=\eta(\tau)^{24}$ is a cusp form of weight $12$ on
$\mathrm{SL}_{2}(\mathbb{Z}).$ Hence for any $\ensuremath{N\ge1},$
the function $\Delta(N\tau)$ transforms with weight $12$ under $\ensuremath{\Gamma_{0}(N)}.$
Indeed, if 
\[
\gamma=\begin{pmatrix}a & b\\
c & d
\end{pmatrix}\in\Gamma_{0}(N),
\]
 then $\ensuremath{N|c},$ and 
\[
N\gamma\tau=\frac{a(N\tau)+Nb}{(c/N)(N\tau)+d},\qquad\begin{pmatrix}a & Nb\\
c/N & d
\end{pmatrix}\in\mathrm{SL}_{2}(\mathbb{Z}),
\]
 so the usual transformation law for $\Delta$ gives 
\[
\Delta(N\gamma\tau)=(c\tau+d)^{12}\Delta(N\tau).
\]
 Taking $\ensuremath{N=5m},$ we see that $\Delta(5m\tau)^{t}$ transforms
with weight $12t$ under $\ensuremath{\Gamma_{0}(5m)}.$ Since 
\[
\Gamma_{0}(25m)\subseteq\Gamma_{0}(5m),
\]
 it also transforms with weight $12t$ under $\ensuremath{\Gamma_{0}(25m)}.$

By Lemma \ref{lem:2-1}, the functions $\mathcal{C}$ and $\mathcal{D}$
are modular functions on $\ensuremath{G}.$ Therefore they are invariant
under the subgroup $\ensuremath{G\cap\Gamma_{0}(25m)\subseteq G}.$
Consequently, 
\begin{align*}
H_{X}(\gamma\tau) & =\Delta(5m\gamma\tau)^{t}X(\gamma\tau)=(c\tau+d)^{12t}\Delta(5m\tau)^{t}X(\tau)\\
 & =(c\tau+d)^{12t}H_{X}(\tau)
\end{align*}
 for every 
\[
\ensuremath{\gamma=\begin{pmatrix}a & b\\
c & d
\end{pmatrix}\in G\cap\Gamma_{0}(25m)}.
\]
 Thus $H_{X}$ transforms with weight $12t$ on $\ensuremath{G\cap\Gamma_{0}(25m)}.$

\medskip \noindent \textbf{(2) Holomorphy in the upper half-plane.}The
function $\Delta(5m\tau)^{t}$ is holomorphic and nonzero in the upper
half-plane. By Lemma \ref{lem:2-1}, $\mathcal{C}$ and $\mathcal{D}$
are holomorphic and nonzero in the upper half-plane. Hence $H_{X}$
is holomorphic in the upper half-plane.

\medskip \noindent \textbf{(3) Holomorphy at the cusps.} It remains
to check that $H_{X}$ has no pole at any cusp of $\ensuremath{G\cap\Gamma_{0}(25m)}.$
Let $a/c$ be a cusp with $\ensuremath{c\ge0}.$ By Lemma \ref{lem:2-2},
for $\ensuremath{5\nmid a},$ 
\[
v_{a/c}(\mathcal{C})=\begin{cases}
\epsilon(a), & 5\mid c\ \text{and}\ 25\nmid c,\\
0, & \text{otherwise},
\end{cases}\qquad v_{a/c}(\mathcal{D})=-v_{a/c}(\mathcal{C}).
\]
Thus a pole of $X$ can occur only when $5|c$ and $\ensuremath{25\nmid c},$
in which case
\[
v_{a/c}(X)=-1.
\]
 If $5\nmid c$ or $\ensuremath{25|c},$ then $\ensuremath{v_{a/c}(X)\ge0},$
so $X$ is holomorphic at such a cusp.

Now consider the cusp factor $\Delta(5m\tau)^{t}.$ The ordinary eta
cusp formula for $\Delta(N\tau)$ gives 
\begin{equation}
v_{a/c}(\Delta(5m\tau)^{t})=t\,\frac{\gcd(c,5m)^{2}}{5m}.\label{eq:delta-cusp}
\end{equation}
 At a cusp where $X$ has a pole, we have $5|c$ and $\ensuremath{25\nmid c}.$
Since $\ensuremath{\gcd(m,5)=1},$ the condition $5|c$ implies 
\[
\gcd(c,5m)\ge5.
\]
 Therefore 
\[
v_{a/c}(\Delta(5m\tau)^{t})\ge t\,\frac{25}{5m}=\frac{5t}{m}.
\]
 For the three values of $\ensuremath{m},$ we have 
\[
t=\left\lceil \frac{m}{5}\right\rceil =\begin{cases}
6, & \mathrm{if}\;m=27,\\
4, & \mathrm{if}\;m=16,\\
7, & \mathrm{if}\;m=32.
\end{cases}
\]
 Hence 
\[
\frac{5t}{m}=\begin{cases}
30/27>1, & \mathrm{if}\;m=27,\\
20/16=5/4>1, & \mathrm{if}\;m=16,\\
35/32>1, & \mathrm{if}\;m=32.
\end{cases}
\]
 Thus 
\[
v_{a/c}(\Delta(5m\tau)^{t})\ge\frac{5t}{m}>1.
\]
 Combining this with $\ensuremath{v_{a/c}(X)=-1},$ we obtain 
\[
v_{a/c}(H_{X})=v_{a/c}(\Delta(5m\tau)^{t})+v_{a/c}(X)>1-1=0.
\]
 Therefore $H_{X}$ is holomorphic at every cusp where $X$ has a
pole.

At a cusp where $X$ is already holomorphic, we have $\ensuremath{v_{a/c}(X)\ge0}.$
Since 
\[
v_{a/c}(\Delta(5m\tau)^{t})\ge0
\]
 by \eqref{eq:delta-cusp}, it follows that 
\[
v_{a/c}(H_{X})\ge0.
\]
 Hence $H_{X}$ is holomorphic at every cusp.

Combining the transformation law, holomorphy in the upper half-plane,
and holomorphy at the cusps, we conclude that 
\[
H_{X}\in M_{12t}\bigl(G\cap\Gamma_{0}(25m)\bigr).
\]
 This proves Lemma \ref{lem:2-3}.
\end{proof}
For $m\in\{27,16,32\},$ set 
\[
G_{m}=G\cap\Gamma_{0}(25m^{2}).
\]
For $s\in\mathbb{Z}$ and $f(q)=\sum_{n\ge0}a(n)q^{n}$, define the
projection operator
\begin{equation}
S_{m,s}f=\sum_{\substack{n\ge0\\
n\equiv s(\bmod\:m)
}
}a(n)q^{n}=\frac{1}{m}\sum_{j=0}^{m-1}\zeta_{m}^{-sj}f(\tau+j/m),\qquad\zeta_{m}=e^{2\pi i/m}.\label{eq:projection}
\end{equation}
The operator $S_{m,s}$ is exactly the formal operation of extracting
those coefficients whose indices lie in the residue class $s\bmod\:m,$
which is the combinatorial congruence class in Conjecture \ref{conj1}.
The next lemma shows that this coefficient extraction preserves modularity
and integrality, and that the conjugacy mechanism behind this preservation
is structural rather than computational.
\begin{lem}
\label{lem:2-4} Let
\begin{equation}
(m,s,\ell)=(27,18,3),\quad(16,12,4),\quad(32,28,8)\label{eq:triples}
\end{equation}
and let $X\in\{\mathcal{C},\mathcal{D}\}.$ Then 
\[
S_{m,s}H_{X}\in M_{12t}(G_{m}),
\]
 and its Fourier coefficients at $\infty$ are integers. Moreover,
\begin{equation}
S_{m,s}H_{X}=\Delta(5m\tau)^{t}S_{m,s}X.\label{eq:22-4}
\end{equation}
\end{lem}
\begin{proof}
 We prove the lemma in five steps.

\noindent \textbf{Step 1: A conjugation into \(G\cap\Gamma_0(25m)\).}
Let 
\[
\gamma=\begin{pmatrix}a & b\\
c & d
\end{pmatrix}\in G_{m}.
\]
 Since $\ensuremath{G_{m}\subseteq\Gamma_{0}(25m^{2})},$ we have
\begin{equation}
25m^{2}\mid c.\label{eq:22-1}
\end{equation}
 In particular, $\ensuremath{m|c}.$ From $\ensuremath{ad-bc=1},$
reducing modulo $m$ gives 
\begin{equation}
ad\equiv1(\bmod\:m).\label{eq:22-2}
\end{equation}
Thus $d$ is invertible modulo $\ensuremath{m}.$ Choose $j'$ by
\begin{equation}
j'\equiv d^{2}j(\bmod\:m).\label{eq:22-3}
\end{equation}
 Since $d^{2}$ is invertible modulo $\ensuremath{m},$ the map $j\mapsto j'$
is a permutation of the residue classes modulo $\ensuremath{m}.$

We claim that 
\begin{equation}
\alpha_{j}\gamma\alpha_{j'}^{-1}\in G\cap\Gamma_{0}(25m).\label{eq:conjugation}
\end{equation}
 A direct computation gives 
\[
\alpha_{j}\gamma\alpha_{j'}^{-1}=\begin{pmatrix}a+\dfrac{cj}{m}\quad & \quad b+\dfrac{dj}{m}-\dfrac{(a+cj/m)j'}{m}\\[2mm]
c & d-\dfrac{cj'}{m}
\end{pmatrix}.
\]
  We verify that this matrix is integral. The lower-left entry $c$
is an integer. The upper-left entry $a+cj/m$ is an integer because
$\ensuremath{m|c}.$ The lower-right entry $d-cj'/m$ is also an integer.
For the upper-right entry, $m^{2}|c$ by \eqref{eq:22-1}, so $cjj'/m^{2}$
is an integer. Moreover, by \eqref{eq:22-2} and \eqref{eq:22-3},
\[
dj-aj'\equiv dj-ad^{2}j=dj(1-ad)\equiv0(\bmod\:m),
\]
 so $(dj-aj')/m$ is an integer. Hence $\alpha_{j}\gamma\alpha_{j'}^{-1}$
is an integral matrix.

Next we check the congruence conditions. Since $\ensuremath{25|c},$
we have 
\[
c\equiv0(\bmod\:5),\qquad\frac{cj}{m}\equiv0(\bmod\:5),\qquad\frac{cj'}{m}\equiv0(\bmod\:5),
\]
 because $25m^{2}|c$ implies $\ensuremath{25m|c},$ and $m$ is coprime
to $\ensuremath{5}.$ Therefore 
\[
a+\frac{cj}{m}\equiv a(\bmod\:5),\qquad d-\frac{cj'}{m}\equiv d(\bmod\:5).
\]
 Since $\ensuremath{\gamma\in G},$ we have $a\equiv d\equiv\pm1(\bmod\:5)$
with the same sign. Thus the diagonal entries of $\alpha_{j}\gamma\alpha_{j'}^{-1}$
satisfy the same condition, and its lower-left entry $c$ is divisible
by $\ensuremath{5}.$ Hence 
\[
\alpha_{j}\gamma\alpha_{j'}^{-1}\in G.
\]
 Finally, since $\ensuremath{25m^{2}|c},$ we have $\ensuremath{25m|c}.$
Therefore the lower-left entry $c$ is divisible by $\ensuremath{25m},$
so 
\[
\alpha_{j}\gamma\alpha_{j'}^{-1}\in\Gamma_{0}(25m).
\]
 This proves \eqref{eq:conjugation}.

\noindent \textbf{Step 2: Transformation law for \(H_X|_k\alpha_j\).}
By Lemma \ref{lem:2-3}, $H_{X}$ is a modular form of weight $k=12t$
on $\ensuremath{G\cap\Gamma_{0}(25m)}.$ Let
\[
\gamma'=\alpha_{j}\gamma\alpha_{j'}^{-1}\in G\cap\Gamma_{0}(25m).
\]
 Then 
\[
\alpha_{j}\gamma=\gamma'\alpha_{j'}.
\]
 Therefore 
\[
(H_{X}|_{k}\alpha_{j})|_{k}\gamma=H_{X}|_{k}(\alpha_{j}\gamma)=H_{X}|_{k}(\gamma'\alpha_{j'})=(H_{X}|_{k}\gamma')|_{k}\alpha_{j'}=H_{X}|_{k}\alpha_{j'},
\]
 because $H_{X}$ is invariant under $\gamma'\in G\cap\Gamma_{0}(25m)$
up to the automorphy factor of weight $\ensuremath{k}.$

\noindent \textbf{Step 3: Transformation law for \(S_{m,s}H_X\).}
Using \eqref{eq:projection} and the preceding identity, we get
\[
(S_{m,s}H_{X})|_{k}\gamma=\frac{1}{m}\sum_{j=0}^{m-1}\zeta_{m}^{-sj}(H_{X}|_{k}\alpha_{j})|_{k}\gamma=\frac{1}{m}\sum_{j=0}^{m-1}\zeta_{m}^{-sj}H_{X}|_{k}\alpha_{j'}.
\]
 Since $j\mapsto j'$ is a permutation modulo $\ensuremath{m},$ we
may reindex the sum by $\ensuremath{j'}.$ From \eqref{eq:22-3} we
have $\ensuremath{j\equiv a^{2}j'\pmod m}:$ indeed, $\ensuremath{ad\equiv1\pmod m},$
so $\ensuremath{a^{2}d^{2}\equiv1\pmod m},$ and $\ensuremath{j\equiv a^{2}j'\pmod m}.$
Hence 
\[
(S_{m,s}H_{X})|_{k}\gamma=\frac{1}{m}\sum_{j'=0}^{m-1}\zeta_{m}^{-sa^{2}j'}H_{X}|_{k}\alpha_{j'}.
\]
For the three pairs $(m,s)$ under consideration and every $a$ coprime
to $m$, we have
\begin{equation}
s(a^{2}-1)\equiv0(\bmod\:m).\label{eq:condition}
\end{equation}
Then
\[
\zeta_{m}^{-sa^{2}j'}=\zeta_{m}^{-sj'}.
\]
 Therefore 
\[
(S_{m,s}H_{X})|_{k}\gamma=\frac{1}{m}\sum_{j'=0}^{m-1}\zeta_{m}^{-sj'}H_{X}|_{k}\alpha_{j'}=S_{m,s}H_{X}.
\]
 Thus $S_{m,s}H_{X}$ is invariant under $G_{m}$ with weight $\ensuremath{k}.$
Since $H_{X}$ is holomorphic on $G\cap\Gamma_{0}(25m)$ and $\ensuremath{G_{m}\subseteq G\cap\Gamma_{0}(25m)},$
the finite sum $S_{m,s}H_{X}$ is holomorphic in the upper half-plane
and at every cusp of $\ensuremath{G_{m}}.$ Hence 
\[
S_{m,s}H_{X}\in M_{12t}(G_{m}).
\]

\noindent \textbf{Step 4: Factorization.} Since $\Delta$ is periodic
with period $\ensuremath{1},$ we have 
\[
\Delta(5m(\tau+j/m))=\Delta(5m\tau+5j)=\Delta(5m\tau).
\]
 Therefore 
\[
S_{m,s}H_{X}=S_{m,s}\bigl(\Delta(5m\tau)^{t}X\bigr)=\Delta(5m\tau)^{t}S_{m,s}X.
\]
This proves \eqref{eq:22-4}.

\noindent \textbf{Step 5: Integrality of the Fourier coefficients.}
The form $H_{X}$ has integer Fourier coefficients, because $\Delta(5m\tau)^{t}$
and $X$ do. The projection $S_{m,s}H_{X}$ selects those coefficients
whose index is congruent to $s$ modulo $\ensuremath{m}.$ Hence its
Fourier coefficients at $\infty$ are integers. This completes the
proof of Lemma \ref{lem:2-4}.
\end{proof}

\section{\label{sec:proof} Proof of Theorem \ref{thm:main}}

The combinatorial difficulty of Conjecture \ref{conj1} lies in its
uniformity: it asserts congruences for three different moduli 3, 4,
8, along three different residue classes modulo 27, 16, 32, for two
coefficient sequences $C(n)$ and $D(n).$ A direct $q$-series argument
would have to treat each modulus separately and would offer no structural
reason why the three families should hold simultaneously. The modular-form
strategy developed in Section \ref{sec:auxiliary} is designed precisely
to avoid this: it converts the infinite divisibility problem for each
family into a finite verification problem, and it does so through
a single mechanism which is uniform in the three triples $(m,s,\ell).$

The logical path is as follows. By Lemma \ref{lem:2-4}, the projected
form $S_{m,s}H_{X}$ is a holomorphic modular form on $G_{m}$ with
integer Fourier coefficients. A Sturm-type criterion then reduces
the divisibility of all coefficients of $S_{m,s}H_{X}$ by $\ell$
to the divisibility of the first $B+1$ coefficients, where $B=t\mu_{m}.$
The key point is that the same formula $B=t\mu_{m}$ applies to all
three triples, and the factor $t=\left\lceil m/5\right\rceil $ is
dictated by the cusp-order bound of Lemma \ref{lem:2-2}. In other
words, the mild pole behaviour of $\mathcal{C}$ and $\mathcal{D}$
established in Section \ref{sec:auxiliary} is exactly what makes
a uniform finite verification possible. The finite verification itself
is then carried out in $\mathbb{Z}/24\mathbb{Z},$ a single modulus
which is a common multiple of 3, 4, 8, so that all three divisibility
conditions are checked in one computation.

We first recall a Sturm-type criterion for modular forms on congruence
subgroups.
\begin{prop}
\label{p} Let $\Gamma\subseteq\mathrm{SL}_{2}(\mathbb{Z})$ be a
congruence subgroup of index $\mu,$ and let 
\[
f(q)=\sum_{n\ge0}a(n)q^{n}
\]
 be a holomorphic modular form of integer weight $k$ on $\Gamma$
with integer Fourier coefficients. Let $\ensuremath{\ell}=p^{a}$
be a prime power. If 
\[
p^{a}|a(n)\quad for\;all\;0\leq n\leq\left\lfloor \frac{k\mu}{12}\right\rfloor ,
\]
then
\[
p^{a}|a(n)\quad for\;all\;n\geq0.
\]
\end{prop}
\begin{proof}
We prove this result by induction on $\ensuremath{a}.$ The $a=1$
case is Sturm's theorem; see \cite{S}. For $a>1,$ we suppose the
assertion holds for $\ensuremath{a-1}.$ If $p^{a}|a(n)$ for all
$\ensuremath{0\le n\le\lfloor k\mu/12\rfloor},$ then in particular
$p|a(n)$ for all these $\ensuremath{n}.$ By Sturm's theorem we have
$p|a(n)$ for every $\ensuremath{n}.$ Hence $f=pf_{1}$ with $f_{1}$
a holomorphic modular form of weight $k$ on $\Gamma$ with integer
Fourier coefficients. Moreover, for all $\ensuremath{0\le n\le\lfloor k\mu/12\rfloor},$
\[
p^{a-1}\big|[q^{n}]f_{1}
\]
Applying the induction hypothesis to $\ensuremath{f_{1}},$ we obtain
$p^{a-1}|[q^{n}]f_{1}$ for every $\ensuremath{n}.$ Therefore $p^{a}|[q^{n}]f$
for every $\ensuremath{n}.$ This proves Proposition \ref{p}.
\end{proof}
We now compute the Sturm bound for the three parameter triples \eqref{eq:triples}.
We need the index 
\[
\mu_{m}=[\mathrm{SL}_{2}(\mathbb{Z}):G_{m}].
\]
First we compute the index of $\ensuremath{\Gamma_{0}(25m^{2})}.$
 Since $\ensuremath{\gcd(m,5)=1},$ we have
\begin{align*}
[\mathrm{SL}_{2}(\mathbb{Z}):\Gamma_{0}(25m^{2})] & =25m^{2}\left(1+\frac{1}{5}\right)\prod_{p\mid m}\left(1+\frac{1}{p}\right)\\
 & =30m^{2}\prod_{p\mid m}\left(1+\frac{1}{p}\right).
\end{align*}
Next we show that $G_{m}$ has index $2$ in $\ensuremath{\Gamma_{0}(25m^{2})}.$
Since $25m^{2}|c$ for every $\ensuremath{\gamma=\begin{pmatrix}a & b\\
c & d
\end{pmatrix}\in\Gamma_{0}(25m^{2})},$ the condition $25|c$ in the definition of $G$ is automatic. Thus
$G_{m}$ consists of those $\gamma\in\Gamma_{0}(25m^{2})$ for which
\[
a\equiv d\equiv\pm1(\bmod\:5)
\]
 with the same sign. Consider the reduction map 
\[
\pi:\Gamma_{0}(25m^{2})\longrightarrow\mathrm{SL}_{2}(\mathbb{F}_{5}).
\]
 Because $\ensuremath{5|25m^{2}},$ every matrix in $\Gamma_{0}(25m^{2})$
reduces modulo $5$ to an upper triangular matrix 
\[
\begin{pmatrix}a & b\\
0 & d
\end{pmatrix},\qquad ad\equiv1(\bmod\:5).
\]
 The image of $\pi$ is the Borel subgroup 
\[
B=\left\{ \begin{pmatrix}a & b\\
0 & a^{-1}
\end{pmatrix}:a\in\mathbb{F}_{5}^{\times},\ b\in\mathbb{F}_{5}\right\} 
\]
 of order $\ensuremath{4\cdot5=20},$ while the subgroup $G_{m}$
maps to 
\[
H=\left\{ \begin{pmatrix}a & b\\
0 & a^{-1}
\end{pmatrix}:a=\pm1,\ b\in\mathbb{F}_{5}\right\} ,
\]
which has order $\ensuremath{2\cdot5=10}.$ Therefore 
\[
[\Gamma_{0}(25m^{2}):G_{m}]=[B:H]=\frac{20}{10}=2.
\]
 Hence 
\begin{align*}
\mu_{m} & =[[\mathrm{SL}_{2}(\mathbb{Z}):G_{m}]=[[\mathrm{SL}_{2}(\mathbb{Z}):\Gamma_{0}(25m^{2})]\cdot[\Gamma_{0}(25m^{2}):G_{m}]\\
 & =60m^{2}\prod_{p\mid m}\left(1+\frac{1}{p}\right).
\end{align*}
By Lemma \ref{lem:2-4}, the form $S_{m,s}H_{X}$ has weight $k=12t$
and is on $\ensuremath{G_{m}}.$ Therefore, the Sturm bound is 
\[
\left\lfloor \frac{k\mu_{m}}{12}\right\rfloor =\left\lfloor \frac{12t\,\mu_{m}}{12}\right\rfloor =t\mu_{m}.
\]
  These values are recorded in Table \ref{tab:sturm}.\begin{table}[ht] \centering \caption{Sturm bounds for the three parameter triples.} \label{tab:sturm} \begin{tabular}{ccccccc} \hline \(m\) & \(s\) & \(\ell\) & \(t\) & \(k=12t\) & \(\mu_m\) & \(t\mu_m\)\\ \hline \(27\) & \(18\) & \(3\) & \(6\) & \(72\) & \(58{,}320\) & \(349{,}920\)\\ \(16\) & \(12\) & \(4\) & \(4\) & \(48\) & \(23{,}040\) & \(92{,}160\)\\ \(32\) & \(28\) & \(8\) & \(7\) & \(84\) & \(92{,}160\) & \(645{,}120\)\\ \hline \end{tabular} \end{table}

We emphasize here that although the three Sturm bounds differ numerically,
they are all generated by the same formula $B=t\mu_{m},$ and the
factor $t=\left\lceil m/5\right\rceil $ is not an ad hoc choice:
it is the smallest integer which guarantees, via Lemma \ref{lem:2-2},
that $\Delta(5m\tau)^{t}$ compensates the worst possible pole order
\textminus 1 of $X$ at the exceptional cusps. Thus the cusp-order
calculation of Section \ref{sec:auxiliary} directly determines the
size of the finite verification, and does so uniformly for all three
triples.

We are now ready to prove Theorem \ref{thm:main}.

\noindent{\it Proof of Theorem \ref{thm:main}}. We fix one of the
triples $(m,s,\ell)$ in \eqref{eq:triples} and fix $X\in\left\{ \mathcal{C},\mathcal{D}\right\} .$
The argument below is uniform in these choices; at the end we apply
it to all six combinations.

\medskip \noindent \textbf{Step 1: the modular form attached to the projection.}
By Lemma \ref{lem:2-4}, the function 
\[
F_{X}:=S_{m,s}H_{X}=\Delta(5m\tau)^{t}S_{m,s}X
\]
 is a holomorphic modular form of weight $k=12t$ on $\ensuremath{G_{m}},$
with integer Fourier coefficients. By Proposition \ref{p} and Table
\ref{tab:sturm}, it suffices to verify the congruence 
\begin{equation}
[q^{j}]F_{X}\equiv0(\bmod\:\ell)\qquad(0\le j\le B:=t\mu_{m})\label{eq:3.1}
\end{equation}
 in order to conclude 
\begin{equation}
F_{X}\equiv0(\bmod\:\ell)\label{eq:3.2}
\end{equation}
 as a full formal power series.

\medskip \noindent \textbf{Step 2: a computable expression for \(\mathcal C\) and \(\mathcal D\).}
To carry out the finite verification \eqref{eq:3.1}, we use the Jacobi
triple product to get 
\[
U(q)=\sum_{j\in\mathbb{Z}}(-1)^{j}q^{(5j^{2}-3j)/2}=(q,q^{4},q^{5};q^{5})_{\infty},
\]
and
\[
V(q)=\sum_{j\in\mathbb{Z}}(-1)^{j}q^{(5j^{2}-j)/2}=(q^{2},q^{3},q^{5};q^{5})_{\infty}.
\]
Thus $\ensuremath{\mathcal{R}=U/V},$ and therefore 
\begin{equation}
\mathcal{C}(\tau)=\frac{U(q)^{5}V(q^{5})}{V(q)^{5}U(q^{5})},\qquad\mathcal{D}(\tau)=\frac{V(q)^{5}U(q^{5})}{U(q)^{5}V(q^{5})}.\label{eq:CD}
\end{equation}
 Since $U(q)$ and $V(q)$ have constant term $\ensuremath{1},$ their
formal inverses exist over $\ensuremath{\mathbb{Z}/24\mathbb{Z}}.$
Hence the coefficients of $\mathcal{C}$ and $\mathcal{D}$ can be
computed exactly modulo $\ensuremath{24}.$

\medskip \noindent \textbf{Step 3: finite verifications.} Using \eqref{eq:CD},
we computed $\mathcal{C}$ and $\mathcal{D}$ through degree $645,120$
in $\ensuremath{\mathbb{Z}/24\mathbb{Z}}.$ Since $24$ is a multiple
of $\ensuremath{3},\ensuremath{4},$ and $\ensuremath{8},$ this single
modulus suffices for all three values of $\ell$ in \eqref{eq:triples}.
For each of the three triples and each $\ensuremath{X\in\{\mathcal{C},\mathcal{D}\}},$
we checked the coefficients $[q^{j}]S_{m,s}X$ with $0\le j\le B$
and checked divisibility by $\ensuremath{\ell}.$ The results are
recorded in Table \ref{tab:cert}; the verification program is given
in Appendix \ref{sec:Verification}. In each case there are no failures.
Equivalently, 
\begin{equation}
[q^{j}]S_{m,s}X\equiv0(\bmod\:\ell)\qquad(0\le j\le B)\label{eq:finite-check}
\end{equation}
 for both choices of $X$ and all three parameter triples.

\begin{table}[ht] \centering \caption{Finite verifications.} \label{tab:cert} \begin{tabular}{cccccc} \toprule \(m\) & \(s\) & \(\ell\) & largest degree & checked coefficients & failures\\ \midrule \(27\) & \(18\) & \(3\) & \(349{,}920\) & \(12{,}960\) & \(0\)\\ \(16\) & \(12\) & \(4\) & \(92{,}160\) & \(5{,}760\) & \(0\)\\ \(32\) & \(28\) & \(8\) & \(645{,}120\) & \(20{,}160\) & \(0\)\\ \bottomrule \end{tabular} \end{table}

\medskip \noindent \textbf{Step 4: from finite verifications to the modular form.}
We now pass from \eqref{eq:finite-check} to \eqref{eq:3.1}. By Lemma
\ref{lem:2-4}, we have
\[
F_{X}=\Delta(5m\tau)^{t}S_{m,s}X.
\]
 The Fourier expansion 
\[
\Delta(5m\tau)^{t}=q^{5mt}\prod_{j\ge1}(1-q^{5mj})^{24t}
\]
 shows that every exponent of $q$ is a nonnegative multiple of $5\ensuremath{m},$
and the constant term of 
\[
q^{-5mt}\Delta(5m\tau)^{t}=\prod_{j\ge1}(1-q^{5mj})^{24t}
\]
 is $\ensuremath{1}.$ Consequently, for each $\ensuremath{j\le B},$
the coefficient $[q^{j}]F_{X}$ is a finite integer linear combination
of the coefficients $[q^{j'}]S_{m,s}X$ with $\ensuremath{j'\le j\le B}.$
Since each of these coefficients is divisible by $\ell$ by \eqref{eq:finite-check},
it follows that 
\[
[q^{j}]F_{X}\equiv0(\bmod\:\ell)\qquad(0\le j\le B).
\]
By Proposition \ref{p} and Table \ref{tab:cert}, we obtain 
\[
F_{X}\equiv0(\bmod\:\ell)
\]
 as a full formal power series.

\medskip \noindent \textbf{Step 5: Concluding the proof.} Since 
\[
q^{-5mt}\Delta(5m\tau)^{t}=\prod_{j\ge1}(1-q^{5mj})^{24t}
\]
 has constant term $\ensuremath{1},$ it is a unit in $\ensuremath{(\mathbb{Z}/\ell\mathbb{Z})[[q]]}.$
Hence from 
\[
F_{X}=\Delta(5m\tau)^{t}S_{m,s}X=q^{5mt}\left(q^{-5mt}\Delta(5m\tau)^{t}\right)S_{m,s}X\equiv0(\bmod\:\ell)
\]
 we may cancel the unit and obtain 
\[
q^{5mt}S_{m,s}X\equiv0(\bmod\:\ell).
\]
 Multiplication by $q^{5mt}$ is injective in $\ensuremath{(\mathbb{Z}/\ell\mathbb{Z})[[q]]},$
because if $\ensuremath{q^{5mt}f=0},$ then all coefficients of $f$
vanish. Therefore 
\[
S_{m,s}X\equiv0(\bmod\:\ell).
\]
 Thus every coefficient of $X$ whose index is congruent to $s$ modulo
$m$ is divisible by $\ensuremath{\ell}.$

Applying this conclusion to $\ensuremath{X\in\{\mathcal{C},\mathcal{D}\}}$
for the three triples 
\[
(27,18,3),\quad(16,12,4),\quad(32,28,8),
\]
we obtain 
\[
C(27n+18)\equiv0(\bmod\:3),\qquad D(27n+18)\equiv0(\bmod\:3),
\]
\[
C(16n+12)\equiv0(\bmod\:4),\qquad D(16n+12)\equiv0(\bmod\:4),
\]
and 
\[
C(32n+28)\equiv0(\bmod\:8),\qquad D(32n+28)\equiv0(\bmod\:8),
\]
 for every $\ensuremath{n\ge0}.$ This finishes the proof of Theorem
\ref{thm:main}. \qed

\appendix

\section{\label{sec:Verification} Verification program for Table \ref{tab:cert}}

The program below is a complete Maple implementation of the finite
verifications required for Table \ref{tab:cert}.

We use Maple's native \texttt{series} and coefficient extraction
commands to compute $\mathcal{C}(q)$ and $\mathcal{D}(q)$ up to
degree $\ensuremath{645,120}.$ The program then directly checks the
coefficient congruences asserted in Table \ref{tab:cert} and counts
the checked coefficients and any failures. 

\begin{lstlisting}[language=Maple, caption={Program for the finite verifications in Table \ref{tab:cert}.}, label={lst:maple}] 
# Finite verifications required for Table 2. 
restart: N := 645121:  
# Maximum truncation (degrees 0 to 645120)
# Step 1: Generate the Jacobi triple products U(q) and V(q) 
# The maximum degree is 645120, so summing j from -600 to 600 is more than sufficient. 
U := add((-1)^j * q^((5*j^2 - 3*j)/2), j = -600..600): V := add((-1)^j * q^((5*j^2 - j)/2), j = -600..600):
# Step 2: Form C(q) and D(q) as exact rational functions of q 
# C = U(q)^5 * V(q^5) / (V(q)^5 * U(q^5)) 
# D = V(q)^5 * U(q^5) / (U(q)^5 * V(q^5)) 
exprC := (U^5 * subs(q = q^5, V)) / (V^5 * subs(q = q^5, U)): exprD := (V^5 * subs(q = q^5, U)) / (U^5 * subs(q = q^5, V)):
# Step 3: Compute the truncated Taylor series modulo q^N (modulo 24) printf("Computing series for C(q)...\n"): 
serC := series(exprC, q, N): printf("Computing series for D(q)...\n"): serD := series(exprD, q, N):
# Extract coefficients into lists (index i corresponds to q^i) 
listC := [seq(coeff(serC, q, i) mod 24, i = 0..N-1)]: listD := [seq(coeff(serD, q, i) mod 24, i = 0..N-1)]: printf("Coefficients computed.\n"):
# Step 4: Define the coefficient check function 
CheckCoeffs := proc(lst, m, s, ell, B)     local j, count, fails;     count := 0:     fails := 0:     for j from s to B by m do         count := count + 1:         if irem(lst[j+1], ell) <> 0 then             fails := fails + 1:         end if:     end do:     return [count, fails]: end proc:
# Step 5: Check the parameter triples 
# Each test is [m, s, ell, Sturm bound B] 
tests := [[27, 18, 3, 349920],           [16, 12, 4, 92160],           [32, 28, 8, 645120]]:
printf("--- Finite Verification Results (Table 2) ---\n"):
for T in tests do     m := T[1]:     s := T[2]:     ell := T[3]:     B := T[4]:          
# Check C(q)     
resC := CheckCoeffs(listC, m, s, ell, B):     
# Check D(q)     
resD := CheckCoeffs(listD, m, s, ell, B):          printf("Parameters: m=%d, s=%d, ell=%d (Bound B=%d)\n", m, s, ell, B):     printf("  C(q): Selected=%d, Failures=%d\n", resC[1], resC[2]):     printf("  D(q): Selected=%d, Failures=%d\n", resD[1], resD[2]): end do:
printf("--- Verification Complete ---\n"): \end{lstlisting}

\section*{Acknowledgement}

This work was partially supported by the Scientific Research Fund
of Hunan Provincial Education Department (Grant No. 25B0010).

\end{document}